\documentclass[12pt]{article}
\usepackage{amsmath,amsthm,amssymb,euscript,verbatim}
\newtheorem{theorem}{Theorem}[section]

\newtheorem{remark}{Remark}[section]

\usepackage{makeidx}         
\usepackage{graphicx}        
\usepackage{multicol}        
\usepackage[bottom]{footmisc}

\makeindex             
\usepackage{latexsym}
\usepackage{epic,eepic,pstricks}

\begin{document}
\title
{\bf  On Cartan's Examples of Isoparametric Hypersurfaces and Their Focal Submanifolds}

\author
{Thomas E. Cecil and Patrick J. Ryan}

\maketitle

\begin{abstract}
This paper is a survey of Cartan's examples of isoparametric hypersurfaces in spheres and their focal submanifolds 
that were described in his fundamental work on the subject, which appeared in four papers
\cite{Car2}--\cite{Car5} published during the period 1938--1940.
\end{abstract}

\section{Isoparametric Hypersurfaces}
\label{sec:isop-hyp-spheres}

We begin with some preliminary notation and definitions, which generally follow those in the book
\cite{CR8}.
By a real {\em space form}
of dimension $n$, we mean a complete, connected, simply connected manifold $\widetilde{M}^n (c)$ with constant sectional curvature $c$.  If $c=0$, then $\widetilde{M}^n (c)$ is the $n$-dimensional Euclidean space ${\bf R}^n$; if
$c=1$, then $\widetilde{M}^n (c)$ is the unit sphere $S^n \subset {\bf R}^{n+1}$; and if $c= -1$, then $\widetilde{M}^n (c)$ is the $n$-dimensional real hyperbolic space $H^n$ (see, for example, \cite[Vol. I, pp. 204--209]{KN}).  For any value of $c>0$, this subject is basically the same as for the case $c=1$, and for any value of $c<0$, it is very similar to the case $c=-1$.  So we will restrict our attention to the cases $c = 0, 1, -1$.

The original definition of a family of isoparametric hypersurfaces in a real space
form $\widetilde{M}^{n+1}$ was formulated in terms of the level sets of an isoparametric function, as we now describe. 
Let $F:\widetilde{M}^{n+1} \rightarrow {\bf R}$ be a non-constant smooth function. The classical Beltrami
differential parameters of $F$ are defined by
\begin{equation}
\label{eq:1-beltrami}
\Delta_1 F = |{\rm grad}\  F|^2, \quad \Delta_2 F = \Delta F\ ({\rm Laplacian}\ {\rm of}\  F),
\end{equation}
where grad $F$ denotes the gradient vector field of $F$.

The function $F$ is said to be {\em isoparametric} if there exist smooth functions $\phi_1$ and $\phi_2$ from ${\bf R}$ to 
${\bf R}$ such that
\begin{equation}
\label{eq:1-beltrami-isoparametric}
\Delta_1 F = \phi_1 \circ F, \quad \Delta_2 F = \phi_2 \circ F.
\end{equation}
That is, both of the Beltrami differential parameters are constant on each level set of $F$. This is the origin of the 
term {\em isoparametric}.
The collection of level sets of an isoparametric function is called an {\em isoparametric family} of sets in
$\widetilde{M}^{n+1}$.

An isoparametric family in ${\bf R}^{n+1}$ consists of either
parallel planes, concentric spheres, or coaxial spherical cylinders,
and their focal sets.
This was first shown for $n=2$ by Somigliana \cite{Som} 
(see also B. Segre \cite{Seg1} and Levi-Civita  \cite{Lev}), and for arbitrary $n$ by B. Segre \cite{Seg}.

In the late 1930's, shortly after the publication of the papers of 
Levi-Civita and Segre,  Cartan \cite{Car2}--\cite{Car5} began a study of isoparametric families in arbitrary
real space forms. 
In Section 2 of his first paper,  Cartan \cite{Car2} showed that an isoparametric family in a real space form $\widetilde{M}$, 
defined as the collection of level sets of an isoparametric function $F$ on $\widetilde{M}$, is locally equivalent to a family of parallel hypersurfaces in $\widetilde{M}$, each of which has constant principal curvatures, together with their focal submanifolds (see \cite[pp. 86--91]{CR8} for the details of this equivalence).

Therefore, in the case where the ambient space is a real space form $\widetilde{M}$,
we will say that  a connected hypersurface $M$ immersed in $\widetilde{M}$ is an {\em isoparametric hypersurface}, if it has constant principal curvatures.

\begin{remark} {\bf Isoparametric hypersurfaces in Riemannian manifolds}
\label{rem-isop-complex-space-forms}
{\rm The definition of an isoparametric hypersurface in the paragraph above is not the appropriate formulation if 
$\widetilde{M}$  is only assumed to be a Riemannian manifold, as can be seen by examples in complex projective space due to Q.-M. Wang \cite{Wang-1a} (see also Thorbergsson \cite{Th6} and \cite[pp. 526--530]{CR8}).
Rather the appropriate generalization is the following: a hypersurface $M$ in a Riemannian manifold $\widetilde{M}$ is {\em isoparametric} if and only if all nearby parallel hypersurfaces $M_t$ have constant mean curvature.

For example, when $\widetilde{M}$ is complex $n$-dimensional projective space ${\bf CP}^n$ or complex
$n$-dimensional hyperbolic space ${\bf CH}^n$, a
(real) isoparametric hypersurface may have non-constant principal curvatures.  Nevertheless, making use of the Hopf map $\pi: S^{2n+1} \rightarrow {\bf CP}^n$ (see \cite[p. 346]{CR8}), one can show that an isoparametric hypersurface in $ {\bf CP}^n$ lifts to an isoparametric hypersurface in $S^{2n+1}$.  This facilitates the classification of isoparametric hypersurfaces in  ${\bf CP}^n$.

In fact, similar methods can be used to classify isoparametric hypersurfaces in ${\bf CH}^n$, although in that case one
must deal with Lorentzian isoparametric hypersufaces in the anti-de Sitter space.
For further details about isoparametric hypersurfaces in complex space forms, see Q.-M. Wang \cite{Wang-1a}, G. Thorbergsson \cite{Th7}, L. Xiao \cite{Xiao}, M. Dominguez-Vazquez et al. \cite{DV2}--\cite{DV3},
and the book \cite[343--531]{CR8}.}
\end{remark}

\section{Cartan's Results}
\label{sec-cartan-formula} 

In Section 3 of his first paper on isoparametric hypersurfaces, Cartan \cite{Car2} derived what he called a
``fundamental formula,'' which is the basis for many of his results on the subject.

The formula involves  the distinct
principal curvatures $\lambda_1,\ldots,\lambda_g$, and their respective multiplicities $m_1,\ldots,m_g$, of an isoparametric
hypersurface $f:M^n \rightarrow \widetilde{M}^{n+1}(c)$ in a space form of constant sectional
curvature $c$.
If the number $g$ of distinct principal curvatures is greater than one, Cartan showed that for each $i$,
$1 \leq i \leq g$, the following equation holds,
\begin{equation}
\label{eq:1-car-id}
\sum_{j\neq i} m_j \frac{c + \lambda_i \lambda_j}{\lambda_i - \lambda_j} = 0.
\end{equation}
Cartan's original formulation of the formula was given in a more complicated form, but he gave the formulation in
equation \eqref{eq:1-car-id} in his third paper \cite[p. 1484]{Car4}.  Equation \eqref{eq:1-car-id} is now known as 
``Cartan's formula.''  Another proof of this formula is given in \cite[pp. 91--96]{CR8}.

Using his formula, Cartan was able to classify isoparametric hypersurfaces in ${\bf R}^{n+1}$ and $H^{n+1}$.  In both cases, one can easily prove using Cartan's formula that the number $g$ of distinct principal curvatures must be either 1 or 2.   Then an elementary argument shows that a connected isoparametric hypersurface is either totally umbilic. i.e., $g=1$, or in the case $g=2$,
it is a an open subset of a tube of constant radius over a totally geodesic submanifold of codimension greater than one (see, for example, \cite[pp. 96--98]{CR8}).

In the sphere $S^{n+1}$, however, Cartan's formula does not lead to the conclusion that $g \leq 2$, and in fact Cartan
produced examples with $g = 1, 2, 3$ and 4 distinct
principal curvatures.  Moreover, he classified isoparametric hypersurfaces $M^n \subset S^{n+1}$ with
$g \leq 3$, as we will discuss in Section \ref{sec-examples-isop-hyp} below.

Cartan \cite[pp. 186--187]{Car2} gave a method for determining the values of the principal curvatures $\lambda_i$ of 
an isoparametric hypersurface $M^n$, if one knows the value of $g$ and the multiplicities $m_1,\ldots,m_g$, as in equation \eqref{eq:1-prin-curv-sph} below,
\begin{equation}
\label{eq:1-prin-curv-sph}
\lambda_i = \cot \theta_i, \ 0 < \theta_i < \pi, \ 1 \leq i \leq g,
\end{equation}
where the $\theta_i$ form an increasing sequence, and $\lambda_i$ has multiplicity $m_i$ on $M$.

Cartan's method is based on 
an analytic calculation involving Cartan's formula \eqref{eq:1-car-id}.
In the case where all of the principal curvatures have the same multiplicity $m = n/g$, Cartan's method 
\cite[pp. 34--35]{Car4} yields the formula \eqref{eq:1-prin-curv-formula-1}
for the principal curvatures given in Theorem 
\ref{thm:1-prin-curv-isop-hyp-cartan} below.

This formula implies that for any point $x \in M$, there
are $2g$ focal points of $(M,x)$ along the normal geodesic to $M$ through $x$, and they are evenly
distributed at intervals of length $\pi/g$.

\begin{theorem} [Cartan]
\label{thm:1-prin-curv-isop-hyp-cartan}  
Let $M \subset S^{n+1}$ be a connected isoparametric hypersurface with $g$ principal curvatures
$\lambda_i = \cot \theta_i$, $0 < \theta_1 < \cdots < \theta_g < \pi$, all having the same multiplicity $m = n/g$. Then
\begin{equation}
\label{eq:1-prin-curv-formula-1}
\theta_i = \theta_1 + (i-1) \frac{\pi}{g} , \quad 1 \leq i \leq g.
\end{equation} 
For any point $x \in M$, there
are $2g$ focal points of $(M,x)$ along the normal geodesic to $M$ through $x$, and they are evenly
distributed at intervals of length $\pi/g$.
\end{theorem}

\begin{remark}
\label{rem-munzner-thm}
{\rm M\"{u}nzner \cite{Mu}
later showed that formula \eqref{eq:1-prin-curv-formula-1} is valid for all isoparametric hypersurfaces in
$S^{n+1}$, even those where the multiplicities are not all the same. M\"{u}nzner's method is different than
the method of Cartan.
M\"{u}nzner's proof is based on the fact
that the set of focal points along a normal geodesic circle to $M \subset S^{n+1}$ is invariant under the 
dihedral group $D_g$
of order $2g$ that acts on the normal circle and is generated by reflections in the focal points.
M\"{u}nzner's proof also yields the result that 
the multiplicities satisfy $m_i = m_{i+2}$ (subscripts mod $g$) (see, for example, \cite[p. 108--109]{CR8}).
Thus there are at most two distinct multiplicities, $m_1, m_2$.}
\end{remark}

Using a lengthy calculation involving equation \eqref{eq:1-prin-curv-formula-1}, Cartan \cite[pp. 364--367]{Car3}
proved the important result
that any isoparametric family with $g$
distinct principal curvatures all having the same multiplicity is algebraic in sense of Theorem \ref{thm:1-Cartan-1} 
below.

Recall that a
function $F:{\bf R}^{n+2} \rightarrow {\bf R}$ is {\em homogeneous of degree} $g$, if 
$g$ is a positive integer such that $F(tx) = t^g F(x)$, for all
$t \in {\bf R}$ and $x \in {\bf R}^{n+2}$.

\begin{theorem} [Cartan]
\label{thm:1-Cartan-1} 
Let $M \subset S^{n+1} \subset {\bf R}^{n+2}$ be a connected isoparametric hypersurface with $g$ principal curvatures
$\lambda_i = \cot \theta_i$, $0 < \theta_1 < \cdots < \theta_g < \pi$, 
all having the same multiplicity $m$.
Then $M$ is an open subset of a level set of the restriction to $S^{n+1}$ of a homogeneous polynomial $F$ 
on ${\bf R}^{n+2}$ of degree $g$ satisfying the differential equations,
\begin{equation}
\label{eq:1-Muenzner-diff-eq-1-car}
\Delta_1 F = |{\mbox{\rm grad }}F|^2 = g^2 r^{2g-2},
\end{equation}
\begin{equation}
\label{eq:1-Muenzner-diff-eq-2-car}
\Delta_2 F = \Delta F = 0 \ (F\ {\mbox{\it is\ harmonic}}),
\end{equation}
where $r = |x|$, for $x \in  {\bf R}^{n+2}$.
\end{theorem}

Cartan's Theorem \ref{thm:1-Cartan-1} was
a forerunner of M\"{u}nzner's general result (see Theorem \ref{thm:1-Muenzner-1} below)
that every isoparametric hypersurface with $g$ distinct principal curvatures
is algebraic (regardless of the multiplicities of the principal curvatures), and its
defining homogeneous polynomial $F$ of degree $g$ satisfies certain conditions on $\Delta_1 F$ and $\Delta_2 F$, 
which generalize the conditions that Cartan found in Theorem \ref{thm:1-Cartan-1}.
(See also the book \cite[pp. 111--130]{CR8} for a discussion of M\"{u}nzner's theorem.)

\begin{theorem}[M\"{u}nzner]
\label{thm:1-Muenzner-1} 
Let $M \subset S^{n+1} \subset {\bf R}^{n+2}$ be a connected isoparametric hypersurface with $g$ 
principal curvatures
$\lambda_i = \cot \theta_i$, $0 < \theta_1 < \cdots < \theta_g < \pi$, with respective multiplicities $m_i$.
Then $M$ is an open subset of a level set of the restriction to $S^{n+1}$ of a homogeneous polynomial $F$ 
on ${\bf R}^{n+2}$ of degree $g$ satisfying the differential equations,
\begin{equation}
\label{eq:1-Muenzner-diff-eq-1}
\Delta_1 F = |{\rm grad}\  F|^2 = g^2 r^{2g-2},
\end{equation}
\begin{equation}
\label{eq:1-Muenzner-diff-eq-2}
\Delta_2 F = \Delta F = c r^{g-2},
\end{equation}
where $r = |x|$, and $c = g^2 (m_2 - m_1)/2$.
\end{theorem}

M\"{u}nzner \cite[p. 65]{Mu}  called $F$ the {\em Cartan polynomial} of $M$.  Now $F$ is usually referred to as the
{\em Cartan-M\"{u}nzner polynomial} of $M$.  Equations (\ref{eq:1-Muenzner-diff-eq-1})--(\ref{eq:1-Muenzner-diff-eq-2})
are called the {\em Cartan-M\"{u}nzner differential equations}.

The following theorem (see \cite[pp. 112--113]{CR8}) relates the gradient and Laplacian of the
restriction $V$ of $F$ to $S^{n+1}$ to the gradient and Laplacian of $F:{\bf R}^{n+2} \rightarrow {\bf R}$ itself.

\begin{theorem}
\label{thm:1-beltrami-homogeneous} 
Let $F:{\bf R}^{n+2} \rightarrow {\bf R}$ be a homogeneous function of degree $g$, and let $V$ be the
restriction of $F$ to $S^{n+1}$.  Then
\begin{enumerate}
\item[${\rm(a)}$] $|{\rm grad}^S V|^2 = |{\rm grad}^E F|^2 - g^2 F^2$,
\item[${\rm(b)}$] $\Delta^S V = \Delta^E F - g (g-1) F - g (n+1) F$,
\end{enumerate}
where ${\rm grad}^S V$ and $\Delta^S V$ denote the gradient and Laplacian of $V$
as a real-valued function on the sphere $S^{n+1}$, and
${\rm grad}^E F$ and $\Delta^E F$ denote the gradient and Laplacian of $F$
as a real-valued function on the space $ {\bf R}^{n+2}$.
\end{theorem}

Using this theorem, one can compute that the
restriction $V$ of $F$ to $S^{n+1}$ satisfies the differential equations,

\begin{equation}
\label{eq:1-Muenzner-diff-eq-S-1}
|{\rm grad}^S V|^2 =  g^2 (1 - V^2),
\end{equation}
\begin{equation}
\label{eq:1-Muenzner-diff-eq-S-2}
\Delta^S V  = c - g(n+g) V,
\end{equation}
where $c = g^2 (m_2 - m_1)/2$.

Thus, if $F$ is a Cartan-M\"{u}nzner polynomial on ${\bf R}^{n+2}$, then its restriction $V$ to $S^{n+1}$
 is an isoparametric function on $S^{n+1}$ in the sense of Cartan 
(as in equation \eqref{eq:1-beltrami-isoparametric}), since both
$|{\rm grad}^S V|^2$ and $\Delta^S V$ are functions of $V$ itself.

Note that Cartan \cite[pp. 364--365]{Car3} proved that if the
principal curvatures all have the same multiplicity $m$, then the polynomial $F$ must be harmonic.  This agrees with 
M\"{u}nzner's condition \eqref{eq:1-Muenzner-diff-eq-2} in the case where all of the multiplicities are equal,
since then the number $c = g^2 (m_2 - m_1)/2$ in equation \eqref{eq:1-Muenzner-diff-eq-2} equals zero,
and so $F$ is harmonic.

\section{Cartan's Examples}
\label{sec-examples-isop-hyp} 
We now begin to describe Cartan's examples of isoparametric hypersurfaces $M \subset S^{n+1} \subset {\bf R}^{n+2}$ 
according to the number $g$ of distinct principal curvatures.  We will describe each example from
various points of view, following the presentations in our books \cite[pp. 294--304]{CR7} and \cite[pp. 144--159]{CR8}
very closely.

Cartan produced examples with $g = 1,2,3$ and $4$ distinct
principal curvatures, and he classified isoparametric hypersurfaces with $g \leq 3$ principal curvatures.  In  his
examples with $g=3$ or $g=4$ principal curvatures, the principal curvatures all have the same multiplicity.

Cartan noted that all of his examples
are homogeneous, each being an orbit of a point under an appropriate closed subgroup of $SO(n+2)$.  
Based on his results and the properties of his examples, Cartan  \cite{Car3} asked the following three questions, all of which were answered in the 1970's, as we will describe below.

\subsubsection*{Cartan's questions}

\begin{enumerate}
\item For each positive integer $g$, does there exist an isoparametric family with $g$ distinct principal
curvatures of the same multiplicity?\\
\item Does there exist an isoparametric family of hypersurfaces with more than three distinct principal curvatures such
that the principal curvatures do not all have the same multiplicity?\\
\item Does every isoparametric family of hypersurfaces admit a transitive group of isometries?
\end{enumerate}

In the early 1970's, Nomizu
\cite{Nom3}--\cite{Nom4} wrote two papers describing the highlights of Cartan's work.  He also generalized 
Cartan's example with four principal curvatures of multiplicity one to produce examples (described 
in Subsection \ref{sub-sec:g=4} below)
with four principal 
curvatures having multiplicities $m_1 = m_3 = m$, and $m_2 = m_4 =1$, for any positive integer $m$.
This answered
Cartan's Question 2 in the affirmative.  

Nomizu also proved that every focal submanifold of every isoparametric
hypersurface is a minimal submanifold of $S^{n+1}$. This also follows from M\"{u}nzner's work \cite{Mu}--\cite{Mu2}, and M\"{u}nzner's proof is different than that of Nomizu.

\begin{remark}
\label{rem-Takagi-Takahashi}
{\rm In 1972, Takagi and Takahashi \cite{TT} gave a complete classification of all homogeneous
isoparametric hypersurfaces
in $S^{n+1}$, based on the work of Hsiang and Lawson \cite{HL}.  Takagi and Takahashi showed that each homogeneous
isoparametric hypersurface in $S^{n+1}$ is a principal orbit of the isotropy representation
of a Riemannian symmetric space of 
rank 2, and they gave a complete list of examples \cite[p.480]{TT}.  This list contains examples with $g=6$
principal curvatures as well as those with $g=1,2,3,4$ principal curvatures. In some cases with $g=4$, the principal
curvatures do not all have the same multiplicity, so this also provided an affirmative answer to 
Cartan's Question 2.}
\end{remark} 

At about the same time as the papers of Nomizu and Takagi-Takahashi, M\"{u}nzner published two preprints that 
greatly extended Cartan's work and have served as the basis for much of the research in the field since that
time.  The preprints were eventually published as papers \cite{Mu}--\cite{Mu2} in 1980--1981. Of course,
one of M\"{u}nzner's primary results is that the number $g$ of distinct principal curvatures of an isoparametric hypersurface in a sphere equals $1,2,3,4$ or 6,
and thus the answer to Cartan's Question 1 is negative.

Finally, the answer to Cartan's Question 3 is also negative, as was first shown by the construction of inhomogeneous
isoparametric hypersurfaces with $g=4$ principal curvatures
by Ozeki and Takeuchi \cite{OT}--\cite{OT2} in 1975.  Their construction was then generalized in 1981 to yield even more inhomogeneous examples by 
Ferus, Karcher and M\"{u}nzner \cite{FKM}.

We now discuss Cartan's examples in detail according to the value of $g$.

\subsection{$g = 1$}
\label{sub-sec:g=1}

If a connected isoparametric hypersurface $M^n \subset S^{n+1} \subset {\bf R}^{n+2}$ has $g=1$ principal curvature, 
then $M^n$ is a totally umbilic hypersurface
in $S^{n+1}$.  As such, $M^n$ is an
open subset of a great or small hypersphere in $S^{n+1}$.  Each great or small hypersphere is the set where $S^{n+1}$
intersects a hyperplane in ${\bf R}^{n+2}$.  The spheres lying in hyperplanes perpendicular to a given diameter of $S^{n+1}$
make up an isoparametric family whose focal set consists of the two end-points (poles) of the diameter.  We now describe these hypersurfaces in terms of the general theory of isoparametric hypersurfaces that we have developed.

Let $p$ be a unit vector in ${\bf R}^{n+2}$ which we will use to determine a diameter of $S^{n+1}$.  Let $F$
be the linear height function,
\begin{equation}
\label{eq:1-height-function}
F(z) = \langle z, p \rangle,
\end{equation}
for $z \in {\bf R}^{n+2}$, where $\langle \ ,\ \rangle$ is the usual inner product on ${\bf R}^{n+2}$.
 Then one easily computes that
\begin{displaymath}
{\rm grad}^E F = p,\quad |{\rm grad}^E F|^2 =1, \quad \Delta^E F = 0,
\end{displaymath}
so that $F$ satisfies the Cartan-M\"{u}nzner differential equations
with $g=1$ and $c=0$, and so $F$ is a Cartan-M\"{u}nzner polynomial.
A calculation using Theorem \ref{thm:1-beltrami-homogeneous}
shows that the restriction $V$ of $F$ to $S^{n+1}$ satisfies
\begin{equation}
\label{eq:1-Beltrami-g-1}
|{\rm grad}^S V|^2 =1 - V^2, \quad \Delta^S V = - (n+1) V,
\end{equation}
so that $V$ is an isoparametric function on $S^{n+1}$.  In fact, it is useful to note that at any 
$z \in S^{n+1}$, we have
\begin{equation}
\label{eq:1-grad-V-1}
{\rm grad}^S V = p - \langle z,p \rangle z.
\end{equation}
We consider the level sets
\begin{displaymath}
M_s = \{z \in S^{n+1} \mid \langle z,p \rangle = \cos s \},\quad 0 \leq s \leq \pi.
\end{displaymath}
Except for the two focal submanifolds,
$M_0 = \{p\}$ and $M_{\pi} = \{-p\}$, each level set is an $n$-sphere with radius $\sin s$ lying in the
hyperplane situated $1 - \cos s$ units below the north pole $p$.  From the point of view of the intrinsic geometry of 
$S^{n+1}$, the set $M_s$ is the geodesic hypersphere $S(p,s)$ centered at the point $p$ having radius $s$.  Note that
$M_s$ may also be regarded as the sphere $S(-p, \pi - s)$.  When $s= \pi/2$, the great sphere $M_s$ is
the unique minimal hypersurface in the isoparametric family.

Finally, the family $M_s$ can also be realized as the set of orbits of a group action as follows.  Let $G_p$ be the subgroup
of $SO(n+2)$ which leaves the pole $p$ fixed.  The group $G_p$ is a naturally embedded copy of $SO(n+1)$, so
that the orbit $G_p z$ has codimension one in $S^{n+1}$ whenever $z \neq \pm p$.  One easily shows that $G_p z$
is a geodesic sphere through $z$ with centers $\pm p$.

\subsection{$g = 2$}
\label{sub-sec:g=2}
Cartan \cite{Car2} showed that an isoparametric hypersurface in the sphere with two distinct
principal curvatures  must be a standard product of two spheres, as we now describe.

We consider  ${\bf R}^{n+2} = {\bf R}^{p+1} \times {\bf R}^{q+1}$, where $p+q = n$.  
For $z = (x,y)$ in ${\bf R}^{p+1} \times {\bf R}^{q+1}$, let 
\begin{equation}
\label{eq:1-Cartan-poly-2}
F(z) = |x|^2 - |y|^2.
\end{equation}
Then we compute
\begin{equation}
\label{eq:1-Beltrami-g-2}
{\rm grad}^E F = 2(x,-y), \quad |{\rm grad}^E F|^2 = 4r^2, \quad \Delta^E F = 2(p-q),
\end{equation}
so that $F$ satisfies the Cartan-M\"{u}nzner differential equations
(\ref{eq:1-Muenzner-diff-eq-1})--(\ref{eq:1-Muenzner-diff-eq-2}) with $g=2$ and $c= 2(p-q)$,
and $F$ is a Cartan-M\"{u}nzner polynomial. 

Note that since $c = g^2 (m_2-m_1)/2 = 2 (m_2-m_1)$, we have $m_2 - m_1 = p - q$.  Since we also
have $m_1 +m_2 = p + q = n$, we conclude that 
\begin{equation}
\label{eq:m-1-m-2}
m_1 = q,\quad m_2 = p.
\end{equation}
The restriction $V$ of $F$ to $S^{n+1}$ satisfies
\begin{eqnarray}
\label{eq:1-Beltrami-g-2-sphere}
{\rm grad}^S V & = & 2 ((1-V)x,-(1+V)y), \quad
|{\rm grad}^S V|^2 = 4 (1 - V^2), \\ 
\Delta^S V & = & 2(p-q) - 2 (n+2)V,\nonumber
\end{eqnarray}
so that $V$ is an isoparametric function on $S^{n+1}$. 

\subsubsection{Level Sets of $F$ in $S^{n+1}$}

Consider the level set of $F$ in $S^{n+1}$,
\begin{displaymath}
M_s = \{z \in S^{n+1} \mid F(z) = \cos 2s \}, \quad 0 \leq s \leq \frac{\pi}{2}.
\end{displaymath}
Each $M_s$ is the Cartesian product of a $p$-sphere of Euclidean radius $\cos s$ and a $q$-sphere of
radius $\sin s$, except for the two focal submanifolds,
\begin{displaymath}
M_0 = \{(x,y) \mid y=0 \} = S^p \times \{0\}, \quad M_{\frac{\pi}{2}} = \{(x,y) \mid x=0 \} = \{0\} \times S^q.
\end{displaymath}
The $M_s$ form a family of parallel hypersurfaces, all of which are tubes over 
each of the focal submanifolds.  Since the focal submanifolds are totally geodesic, it is easy to compute the shape
operator of $M_s$ using the formula for the shape operator of a tube (see, for example,
\cite[pp. 17--18]{CR8}).  
By considering $M_s$ as a tube of radius $s$ over the totally geodesic submanifold $S^p \times \{0\}$, we get that $M_s$ has two constant principal curvatures,
\begin{displaymath}
\lambda = \cot \left(\frac{\pi}{2} -s\right) = \tan s, \quad \mu = - \cot s,
\end{displaymath}
with respective multiplicities $p$, $q$.

The family $\{M_s\}$ can also be realized as the set of orbits of a group action.  Considering the decomposition
${\bf R}^{n+2} = {\bf R}^{p+1} \times {\bf R}^{q+1}$, the group 
\begin{displaymath}
G = SO(p+1) \times SO(q+1)
\end{displaymath}
is naturally embedded as a subgroup of $SO(n+2)$. Take
any point $z = (x,y)$ with both $x$ and $y$ non-zero.  Then the isotropy subgroup of the $G$-action at $z$ is
isomorphic to $SO(p) \times SO(q)$, so that the orbit $M$ of $z$ has codimension 2 in ${\bf R}^{n+2}$, and hence
codimension 1 in $S^{n+1}$.  

Considering the action of $G$ on each factor, it is clear that the orbit consists
of the product of standard spheres of radii $|x|$ and $|y|$, and hence it coincides with $M_s$, where 
$\cos s = |x|$ and $\sin s = |y|$.  Furthermore, the orbit of $(x,0)$ is the focal submanifold $M_0$,
and the orbit of $(0,y)$ is the focal submanifold $M_{\pi/2}$.

\subsection{$g = 3$}
\label{sub-sec:g=3}
In his paper \cite{Car3}, Cartan gives a complete classification 
of isoparametric hypersurfaces in $S^{n+1}$ with $g=3$ principal curvatures,
and he describes each example in great detail.

Let $M^n \subset S^{n+1} \subset {\bf R}^{n+2}$ be a connected, oriented hypersurface  with field of unit normals $\xi$
having $g=3$ distinct principal curvatures at every point.  The first step in Cartan's classification theorem is to prove the following theorem \cite[pp. 359--360]{Car3}.

\begin{theorem}
\label{thm:isop-hyp-g-3-multiplicities} 
Let $M^n \subset S^{n+1}$  be a connected isoparametric hypersurface with $g=3$ principal curvatures at each point.  Then all of the principal curvatures have the same multiplicity $m$.
\end{theorem}

Cartan proves this using the method of moving frames and conditions on the structural formulas that follow from the isoparametric condition and the assumption that $g=3$, and we refer the reader to Cartan's paper
\cite[pp. 359--360]{Car3} for the proof.

Next, since all of the multiplicities are equal, 
Cartan's Theorem \ref{thm:1-Cartan-1} mentioned
above implies that $M^n$ is an open subset of a level set of the restriction to $S^{n+1}$ of a harmonic
homogeneous polynomial $F$ 
on ${\bf R}^{n+2}$ of degree $g=3$ satisfying the differential equations stated in Theorem \ref{thm:1-Cartan-1}.

Cartan's approach to classifying isoparametric hypersurfaces with
$g=3$ 
is to try to directly determine all homogeneous harmonic polynomials $F$ of degree $g=3$ that satisfy the differential equation from Theorem \ref{thm:1-Cartan-1},

\begin{equation}
\label{eq:1-Muenzner-diff-eq-1-a}
\Delta_1 F = |{\mbox{\rm grad }}F|^2 = g^2 r^{2g-2} = g^2 (x_1^2 + x_2^2 + \cdots + x_{n+2}^2)^{{g-1}},
\end{equation}
where $x = (x_1, x_2, \ldots, x_{n+2}) \in {\bf R}^{n+2}$. 
After a long calculation, Cartan shows that the algebraic determination of such a harmonic polynomial requires the possibility of solving
the algebraic problem:\\

\noindent
{\bf Problem 1:} {\it Represent the product
\begin{displaymath}
(u_1^2 + u_2^2 + \cdots + u_m^2) (v_1^2 + v_2^2 + \cdots + v_m^2)
\end{displaymath}
of two sums of $m$ squares by the sum of the squares of $m$ bilinear combinations of the $u_i$ and the $v_j$}.\\

\noindent
By a theorem of Hurwitz \cite{Hurwitz} on normed linear algebras, this is only possible if $m = 1, 2, 4,$ or 8.

Cartan \cite[pp. 340--341]{Car3} relates the solution of Problem 1 above to the theory of Riemannian spaces admitting an isogonal absolute parallelism, and he uses results from a joint paper of himself and J. A. Schouten 
\cite{cartan-schouten} on that topic to determine the possibilities for $F$.

Ultimately, Cartan \cite{Car3} shows that the required harmonic homogeneous polynomial $F$ of degree 3 
on ${\bf R}^{3m+2}$ must be of the form
$F(x,y,X,Y,Z)$  given by

\begin{equation}
\label{eq:1-Cartan-poly-g-3}
x^3 - 3xy^2 + \frac{3}{2} x(X \overline{X} + Y \overline{Y} - 2 Z \overline{Z}) + \frac{3 \sqrt{3}}{2}
y (X\overline{X} - Y \overline{Y}) + \frac{3 \sqrt{3}}{2} (XYZ + \overline{Z} \overline{Y} \overline{X}).
\end{equation}
In this formula, $x$ and $y$ are real parameters, while $X,Y,Z$ are coordinates in the division algebra
${\bf F} = {\bf R}, {\bf C}, {\bf H}$ (quaternions), ${\bf O}$ (Cayley numbers),
for $m = 1,2,4,8$,
respectively.  

Note that the sum $XYZ + \overline{Z} \overline{Y} \overline{X}$ is twice the real part
of the product $XYZ$.  In the case of the Cayley numbers, multiplication is not associative, but the
real part of $XYZ$ is the same whether one interprets the product as $(XY)Z$ or $X(YZ)$. 

The isoparametric
hypersurfaces in the family are the level sets $M_t$ in $S^{3m+1}$ determined by the equation
$F = \cos 3t$, $0 < t < \pi/3$, where $F$ is the polynomial in
equation \eqref{eq:1-Cartan-poly-g-3}.  The focal submanifolds are obtained by taking $t=0$ and $t= \pi/3$.  These 
focal submanifolds are
a pair of antipodal standard embeddings of the projective plane
${\bf FP}^2$, for the appropriate division
algebra ${\bf F}$ listed above
(see \cite[pp. 74--78]{CR8} for 
more detail on the standard embeddings of projective planes).

For the cases
${\bf F} = {\bf R}, {\bf C}, {\bf H}$, Cartan gave a specific parametrization of the focal submanifold $M_0$
defined by the condition $F=1$ as follows:
\begin{eqnarray}
\label{eq:standard-embedding-focal-set}
X = \sqrt{3} v \overline{w}, \quad Y = \sqrt{3} w \overline{u}, \quad Z = \sqrt{3} u \overline{v},\\
x = \frac{\sqrt{3}}{2}(|u|^2 - |v|^2), \quad y = |w|^2 - \frac{|u|^2 + |v|^2}{2},\nonumber
\end{eqnarray}
where $u,v,w$ are in ${\bf F}$, and $|u|^2 + |v|^2 + |w|^2 = 1$.  This map is invariant under the
equivalence relation
\begin{displaymath}
(u,v,w) \sim (u\lambda,v\lambda,w\lambda), \quad \lambda \in {\bf F},\  |\lambda| = 1.
\end{displaymath}
Thus, it is well-defined on ${\bf FP}^2$, and it is easily shown to be injective on ${\bf FP}^2$. Therefore, it is
an embedding of ${\bf FP}^2$ into $S^{3m+1}$. (This parametrization differs slightly from that given in
\cite[pp. 74--78]{CR8} for the standard embeddings of projective spaces.)

 Cartan states that he does not know of an analogous representation of the focal submanifold in the case $m = 8$. 
 In that case, if  $u, v, w$ are taken to be Cayley numbers, then
the formulas \eqref{eq:standard-embedding-focal-set} are not well defined on the Cayley projective plane.
For example, the product $v \overline{w}$ is not preserved, in general, if we replace $v$ by  $v\lambda$, and $w$ by $w\lambda$, where $\lambda$ is a unit Cayley number.  

Note that even though there is no parametrization corresponding to \eqref{eq:standard-embedding-focal-set},
the Cayley projective plane can be described (see, for example,  Kuiper \cite{Ku3}, Freudenthal \cite{Freu})
as the submanifold
\begin{displaymath}
V = \{A \in M_{3 \times 3} ({\bf O}) \mid \overline{A}^T = A = A^2,\ {\rm trace}\ A = 1 \},
\end{displaymath}
where $M_{3 \times 3} ({\bf O})$ is the space of $3 \times 3$ matrices of Cayley numbers.  This submanifold
$V$ lies in a sphere $S^{25}$ in a $26$-dimensional real subspace of $M_{3 \times 3} ({\bf O})$.

In his paper, Cartan \cite[pp. 342--358]{Car3} writes a separate section for each of the cases $m = 1,2,4,8$, that is, ${\bf F} = {\bf R}, 
{\bf C}, {\bf H},{\bf O}$.  He gives extensive details in each case, describing many remarkable properties, especially in the case $m=8$. 
In particular, he shows that in each case the isoparametric hypersurfaces and the two focal
submanifolds are homogeneous, that is, they are orbits of points in $S^{n+1}$ under the action
of a closed subgroup of $SO(n+2)$.

Cartan's results show that up to congruence,
there is only one isoparametric family of hypersurfaces with $g=3$ principal curvatures for each value of 
$m = 1,2,4,8$.  For each value of $m$, the two focal submanifolds are a pair of antipodal standard embeddings
of the projective plane ${\bf FP}^2$, for the division algebra
${\bf F} = {\bf R}, {\bf C}, {\bf H}$ (quaternions), ${\bf O}$ (Cayley numbers),
for $m = 1,2,4,8$, respectively.  Each isoparametric hypersurface in the family is a tube of constant radius
over each of the two focal submanifolds.

Cartan's classification is closely related to various characterizations of the standard embeddings of ${\bf FP}^2$.
(See Ewert \cite{Ewert}, Little \cite{Little}, and Knarr-Kramer \cite{K-K}.)
For alternative proofs of Cartan's 
classification of isoparametric hypersurfaces with $g=3$ principal curvatures,
see the papers of Knarr and Kramer \cite{K-K}, and Console and Olmos \cite{Console-Olmos-98}.
In a related paper, Sanchez \cite{Sanchez} studied Cartan's isoparametric hypersurfaces from an algebraic point of view. 
(See also the paper of Giunta and Sanchez \cite{Giunta-Sanchez-2014}.)

\subsubsection{Homogeneity of Cartan's hypersurfaces with $g=3$} 

Cartan showed that in each case $m = 1,2,4$ or 8, the isoparametric hypersurfaces and the two focal
submanifolds are homogeneous, that is, they are orbits of points in $S^{n+1}$ under the action
of a closed subgroup of $SO(n+2)$.  Here we give a presentation of this fact for the case ${\bf F} = {\bf R}$,
as in \cite[pp. 297--299]{CR7} or \cite[pp. 154--155]{CR8}.
An analogous construction can be made for the other algebras.  

We consider ${\bf R}^9$ as the space
$M_{3 \times 3} ({\bf R})$ of $3 \times 3$ real matrices with standard inner product
\begin{displaymath}
\langle A,B \rangle = {\rm trace}\ AB^T,
\end{displaymath}
where $B^T$ denotes the transpose of $B$.
We consider the 5-dimensional subspace ${\bf R}^5$ of symmetric matrices with trace zero, and let $S^4$
be the unit sphere in ${\bf R}^5$.  That is,
\begin{displaymath}
S^4 = \{A \in M_{3 \times 3} ({\bf R}) \mid A = A^T,\ {\rm trace}\ A = 0,\ |A| = 1\}.
\end{displaymath}
The group $SO(3)$ acts on $S^4$ by conjugation.  This action is isometric and thus preserves $S^4$.  For every
$A \in S^4$, there exists a matrix $U \in SO(3)$ such that $UAU^T$ is diagonal.  In fact, a direct
calculation shows that every orbit of this action contains a representative of the form $B_t$, where
$B_t$ is a diagonal matrix whose diagonal entries are
\begin{displaymath}
\sqrt{\frac{2}{3}}\  \{\cos (t - \frac{\pi}{3}),\ \cos (t + \frac{\pi}{3}),\ \cos (t + \pi) \}.
\end{displaymath}
If all the eigenvalues of $B_t$ are distinct, then the orbit of $B_t$ is 3-dimensional.  For example, consider
\begin{displaymath}
B_{\pi/6} = {\rm diagonal}\ \{ 1/\sqrt{2}, 0, -1/\sqrt{2}\}.
\end{displaymath}
The isotropy subgroup of $B_{\pi/6}$ under this group action is the set of matrices in $SO(3)$ that commute
with $B_{\pi/6}$.  One can easily compute that this group consists of diagonal matrices in $SO(3)$ with entries
$\pm1$ along the diagonal.  This group is isomorphic to the Klein 4-group ${\bf Z}_2 \times {\bf Z}_2$, and thus
the orbit $M_{\pi/6}$ is isomorphic to $SO(3)/{\bf Z}_2 \times {\bf Z}_2$.  The hypersurface $M_{\pi/6}$ is the unique
minimal hypersurface in the isoparametric family.

The two focal submanifolds are lower-dimensional orbits, and they occur when $B_t$ has a repeated eigenvalue.  For
example, when $t=0$, we have
\begin{displaymath}
B_0 = {\rm diagonal}\ \{ 1/\sqrt{6}, 1/\sqrt{6}, -2/\sqrt{6}\}.
\end{displaymath}
The isotropy subgroup for $B_0$ is the subgroup $S(O(2) \times O(1))$ consisting of all matrices of the form
\begin{displaymath}
\left[ \begin{array}{cc}A&0\\0&\pm 1\end{array}\right], \quad A \in O(2), 
\end{displaymath}
having determinant one.  Thus, $M_0$ is diffeomorphic to $SO(3)/S(O(2) \times O(1))$ which is the real
projective plane ${\bf RP}^2$.  In fact, $M_0$ is a standard embedding of ${\bf RP}^2$, as noted above.  The other focal
submanifold $M_{\pi/3}$ is the orbit of 
\begin{displaymath}
B_{\pi/3} = {\rm diagonal}\ \{ 2/\sqrt{6}, -1/\sqrt{6}, -1/\sqrt{6}\}.
\end{displaymath}
This is also a standard embedding of ${\bf RP}^2$antipodal to $M_0$.

\subsection{$g = 4$}
\label{sub-sec:g=4}

In his fourth paper on the subject, Cartan \cite{Car5} continued
his study of isoparametric hypersurfaces with the property that all the principal curvatures have the same multiplicity $m$.
In that paper, Cartan produced a family of 
isoparametric hypersurfaces with $g=4$ principal curvatures of 
multiplicity $m=1$ in $S^5$,
and a family with $g=4$ principal curvatures of multiplicity $m=2$ in $S^9$.  

In both cases, Cartan describes the geometry and topology of the hypersurfaces and the two focal submanifolds in great detail, pointing out several notable properties.  
(See also the papers of Nomizu \cite{Nom3}--\cite{Nom4} and the books \cite[pp. 299--303]{CR7},
\cite[pp. 155--159]{CR8}, for more detail.)  Cartan also proved that the hypersurfaces and the focal submanifolds are homogeneous in both cases, $m=1$ and $m=2$.

Cartan stated that he could show by a very long calculation (which he did not give) that for an isoparametric hypersurface
with $g=4$ principal curvatures having the same multiplicity $m$, the only possibilities are $m=1$ and $m=2$.  This was shown to be true
later by Grove and Halperin \cite{GH}.
Furthermore, up to congruence, Cartan's examples
given in  \cite{Car5} are the only possibilities for  $m=1$ and $m=2$.  This was shown to be true 
by Takagi \cite{Takagi} in the case $m=1$, and by Ozeki and Takeuchi \cite{OT}--\cite{OT2} in the case $m=2$.  

As noted in Remark \ref{rem-Takagi-Takahashi}, 
in 1972 R. Takagi and T. Takahashi \cite{TT}
published a complete classification of homogeneous isoparametric hypersurfaces in spheres based on
the work of Hsiang and Lawson \cite{HL}.  Takagi and Takahashi showed that each homogeneous isoparametric
hypersurface $M$ in $S^{n+1}$ is a principal orbit of the isotropy representation of a Riemannian symmetric 
space of rank 2, and they gave a complete list of examples \cite[p. 480]{TT}, including
examples with $g$ equal to each of the values $1,2,3,4,$ and 6.

\subsubsection{A homogeneous example with $g=4$}

We begin our treatment of Cartan's isoparametric hypersurfaces with $g=4$ principal curvatures with  class of
homogeneous examples.
This construction is due in full generality to Nomizu \cite{Nom3}--\cite{Nom4}, and it was first
given by Cartan \cite{Car5} in the case of $g=4$ principal curvatures of multiplicity 1.
Our treatment of this example follows Nomizu \cite{Nom3}--\cite{Nom4}, and we follow the text in our
book \cite[pp. 155--159]{CR8} very closely.

These examples are noteworthy, because they are the simplest cases of a large
class of isoparametric hypersurfaces, whose construction is based on representations of Clifford algebras. 
That construction was first given in certain cases by Ozeki and Takeuchi \cite{OT}--\cite{OT2}, 
and then later in full generality by
Ferus, Karcher and M\"{u}nzner \cite{FKM}.  Thus, we will say that such isoparametric hypersurfaces are of 
{\em OT-FKM type}.

We consider the $(m+1)$-dimensional complex vector space ${\bf C}^{m+1}$ as a real vector space
${\bf C}^{m+1} = {\bf R}^{m+1} \oplus i {\bf R}^{m+1}$.  The real inner product on ${\bf C}^{m+1}$ is given by
\begin{displaymath}
\langle z,w \rangle = \langle x,u \rangle + \langle y,v \rangle,
\end{displaymath}
for $z = x + iy$, $w = u + i v$ for $x,y,u,v \in {\bf R}^{m+1}$.  The unit sphere in ${\bf C}^{m+1}$ is 
\begin{displaymath}
S^{2m+1} = \{z \in {\bf C}^{m+1} \mid |z| = 1 \}.
\end{displaymath}
In the following construction, we assume that $m \geq 2$.  In the case $m=1$, this construction reduces to
a product of two circles in $S^3$.
Consider the homogeneous polynomial $F$ of degree 4 on ${\bf C}^{m+1}$  given by
\begin{equation}
\label{eq:1-F-Nomizu}
F(z) = | \sum_{k=0}^m z_k^2 |^2 = (|x|^2 - |y|^2)^2 + 4\  \langle x,y \rangle^2, \ {\rm for}\ z = x + iy.
\end{equation}
A direct computation shows that
\begin{equation}
\label{eq:1-Beltrami-g-4}
|{\rm grad}^E F|^2 = 16 r^2 F, \quad \Delta^E F = 16 r^2,
\end{equation}
and therefore that the restriction $V$ of $F$ to $S^{2m+1}$ satisfies
\begin{equation}
\label{eq:1-Beltrami-g-4-sphere}
|{\rm grad}^S V|^2 = 16 V (1 - V), \quad \Delta^S V = 16 - V (16 + 8m).
\end{equation}
Thus, $V$ is an isoparametric function on $S^{2m+1}$.

\begin{remark}
\label{rem:1-Cartan-poly-Nom}
{\rm Note that $F$ does not satisfy the Cartan-M\"{u}nzner differential equations 
(\ref{eq:1-Muenzner-diff-eq-1})--(\ref{eq:1-Muenzner-diff-eq-2}).  However, as noted by Takagi \cite{Takagi},
the polynomial $\widetilde{F} = r^4 - 2F$ has the same level sets as $F$ on $S^{2m+1}$, since the restriction
$\widetilde{V}$ of $\widetilde{F}$ to $S^{2m+1}$ satisfies $\widetilde{V} = 1 - 2V$.  The function $\widetilde{F}$ satisfies
the equations,
\begin{equation}
\label{eq:1-Beltrami-g-4-tilde}
|{\rm grad}^E \widetilde{F}|^2 = 16 r^6, \quad \Delta^E \widetilde{F} = 8 (m-2) r^2,
\end{equation}
and so it does satisfy the Cartan-M\"{u}nzner differential equations 
with $g=4$ and multiplicities $m_2 = m-1$ and $m_1 = 1$, since $c = g^2 (m_2 - m_1)/2 = 8 (m-2)$.}
\end{remark}

\subsubsection{Focal submanifolds of this family}

We now continue our discussion of this example using the functions $F$ and $V$.
From equation (\ref{eq:1-Beltrami-g-4-sphere}), we see that the focal submanifolds occur when $V=0$ or 1. 
From equation (\ref{eq:1-F-Nomizu}), we that $V=1$ is equivalent to the condition,
\begin{displaymath}
| \sum_{k=0}^m z_k^2 | = 1.
\end{displaymath}
This is easily seen to be equivalent to the condition that $z$ lies in the set
\begin{displaymath}
M_0 = \{e^{i\theta} x \mid x \in S^m \},
\end{displaymath}
where $S^m$ is the unit sphere in the first factor ${\bf R}^{m+1}$.  For $x \in S^m$, we have
\begin{displaymath}
T_x M_0 = T_x S^m \oplus {\rm Span}\   \{ix\}.
\end{displaymath}
Thus, the normal space to $M_0$ at $x$ is,
\begin{displaymath}
T_x^\perp M_0 = \{iy \mid y \in S^m,\ \langle x,y \rangle = 0 \}.
\end{displaymath}
The normal geodesic to $M_0$ through $x$ in the direction $iy$ can be parametrized as
\begin{equation}
\label{eq:1-7.3}
\cos t \ x + \sin t \  iy.
\end{equation}
At the point $e^{i \theta} x$ in the focal submanifold $M_0$, one can easily show that 
\begin{displaymath}
T^\perp M_0 = \{e^{i \theta} y \mid y \in S^m,\ \langle x,y \rangle = 0 \}.
\end{displaymath}
Thus, the normal geodesic to $M_0$ through the point $e^{i \theta} x$ in the direction $e^{i \theta} iy$ can
be parametrized as,
\begin{equation}
\label{eq:1-7.4}
\cos t \ e^{i \theta} x + \sin t \  e^{i \theta} iy = e^{i \theta} (\cos t \ x + \sin t \  iy).
\end{equation}
Let $V_{m+1,2}$ be the Stiefel manifold of orthonormal pairs of vectors $(x,y)$ in ${\bf R}^{m+1}$.
By equations (\ref{eq:1-7.3}) and (\ref{eq:1-7.4}), we see that the tube $M_t$ of radius $t$ over the focal submanifold 
$M_0$ is given by
\begin{equation}
\label{eq:1-tube-Nomizu}
M_t = \{e^{i \theta} (\cos t \ x + \sin t \  iy) \mid (x,y) \in V_{m+1,2} \}.
\end{equation}
In fact, the map $f_t:S^1 \times V_{m+1,2} \rightarrow S^{2m+1}$ given by
\begin{equation}
\label{eq:1-7.5}
f_t (e^{i \theta}, (x,y)) = e^{i \theta} (\cos t \ x + \sin t \  iy),
\end{equation}
is an immersion that is a double covering of the tube $M_t$, since 
\begin{displaymath}
f_t (e^{i \theta}, (x,y)) = f_t (e^{i (\theta + \pi)}, (-x,-y)).
\end{displaymath}
Substituting equation (\ref{eq:1-7.5}) into the defining formula (\ref{eq:1-F-Nomizu}) for $F$ shows that for
$z \in M_t$, the restriction $V$ of $F$ to $S^{2m+1}$ satisfies,
\begin{displaymath}
V(z) = (\cos^2 t - \sin^2 t)^2 = \cos^2 2t.
\end{displaymath}
Therefore, the other focal submanifold determined by the equation $V=0$ occurs when $t = \pi/4$.
By equations (\ref{eq:1-tube-Nomizu}), we get that the focal submanifold $M_{\pi/4}$ consists of points of the form 
$e^{i \theta} (x+iy)/\sqrt{2}$ for $(x,y) \in V_{m+1,2}$.  On the other hand, the equation $V=0$ implies that
$|x| = |y|$ and $\langle x, y \rangle = 0$, and we conclude that
\begin{displaymath}
M_{\pi/4} = \{(x+iy)/\sqrt{2} \mid (x,y) \in V_{m+1,2} \}.
\end{displaymath}
This is an embedded image of the Stiefel manifold $V_{m+1,2}$ which has dimension $2m-1$.  

\begin{remark}
\label{rem:Stiefel-manifold}
{\rm The fact that one of the focal
submanifolds is a Stiefel manifold illustrates an important feature of the general construction of 
Ferus, Karcher and M\"{u}nzner
\cite{FKM}. In that construction, one of the focal submanifolds is always a so-called Clifford-Stiefel manifold
determined by a corresponding Clifford algebra (see also Pinkall-Thorbergsson \cite{PT1}, and the book
\cite[pp. 162--180]{CR8} for more detail.)}
\end{remark}

Since adjacent focal points along a normal geodesic are at a distance $\pi/4$ apart, we know from
M\"{u}nzner's Theorem on the values of the principal curvatures of an isoparametric hypersurface
\cite[p. 108]{CR8} that the tube 
 $M_t$ of radius $t$ over the focal submanifold  $M_0$ 
is an
isoparametric hypersurface with four principal curvatures,
\begin{equation}
\label{eq:1-7.6}
\cot t,\ \cot \left(t+\frac{\pi}{4}\right),\ \cot \left(t+\frac{\pi}{2}\right),\ \cot \left(t+\frac{3\pi}{4}\right).
\end{equation}
The focal submanifolds $M_0$ and $M_{\pi/4}$ have respective dimensions $m+1$ and $2m-1$. Thus, the principal
curvatures in equation (\ref{eq:1-7.6}) have respective multiplicities $(m-1),1,(m-1),1$.  This agrees with
the information that we obtained from the Cartan-M\"{u}nzner polynomial $\widetilde{F}$ 
in Remark \ref{rem:1-Cartan-poly-Nom}.

From equation (\ref{eq:1-7.5}), we see that $M_t$ admits a transitive group of isometries isomorphic
to $SO(2) \times SO(m+1)$, and hence each $M_t$ is an orbit hypersurface.
This is the fifth example on
the list of Takagi and Takahashi \cite{TT} of homogeneous isoparametric hypersurfaces.
Later Takagi \cite{Takagi} showed that if an isoparametric hypersurface $M$ in $S^{2m+1}$ has four principal curvatures with multiplicities $(m-1),1,(m-1),1$, then $M$ is
congruent to a hypersurface $M_t$ in this example.  In particular, $M$ is homogeneous.

Regarding the other focal submanifold $M_0$, the map $f_0:S^1 \times S^m \rightarrow M_0$ given by
$f_0 (e^{i \theta},x) = e^{i \theta}x$ is a double covering of $M_0$, since
\begin{displaymath} 
f_0 (e^{i \theta},x) = f_0 (e^{i (\theta + \pi)},-x).
\end{displaymath}
Hence, $M_0$ can be considered as a quotient 
manifold with identifications given by the map $f_0$.  The two spheres $\{1\} \times S^m$ and
$\{-1\} \times S^m$ are attached via the antipodal map of $S^m$.  Thus, $M_0$ is orientable if and only if the 
antipodal map on $S^m$ preserves orientation on $S^m$.  Therefore, $M_0$ is orientable if $m$ is odd and non-orientable 
if $m$ is even.  

This example illustrates M\"{u}nzner's discussion of the orientability of the focal submanifolds
in Theorem C \cite[p. 59]{Mu}.  M\"{u}nzner shows that in the case $g=4$, there are two possibilities.
In the first situation, one of the focal submanifolds $M_+, M_-$ is orientable and the other is not, in which
case, if $M_+$ is orientable, then the multiplicity $m_+ = 1$.  In the second situation,
both focal submanifolds are orientable, in which case $m_+ + m_-$ is odd or else $m_+ = m_-$ is even.

In our example, $m_+ =1$, $m_- = m-1$ and $m_+ + m_- = m$.  Thus, if $m$ is even, M\"{u}nzner's
theorem implies that the focal submanifold $M_- = M_0$ corresponding to the principal curvature of multiplicity
$m-1$ is non-orientable, as we have shown above.  In the case where $m$ is odd, then both focal submanifolds
in our example are orientable, which is consistent with M\"{u}nzner's theorem, although not a consequence of it.

\begin{remark}
\label{rem:(2,2)-case}
{\rm We will not describe Cartan's example having $g=4$ principal curvatures of multiplicity $m=2$ in detail,
and we refer the reader to Cartan's original paper \cite{Car5}, in which that example is
described in full detail.  Cartan's example is also noteworthy, because it is not of OT-FKM type.  It is one of only
two examples with $g=4$ that is not of OT-FKM type. The other one is the homogeneous example with $g=4$ principal curvatures, and multiplicities
$(m_1, m_2) = (4,5)$.} 
\end{remark}

\section{Classification Results}
\label{chap:2}

Due to the work of many mathematicians, isoparametric hypersurfaces in $S^{n+1}$ have been completely classified, and we will briefly state the classification results now.  As noted above, M\"{u}nzner
proved that the number $g$ of distinct principal curvatures must be $1,2,3, 4$ or 6.
We now summarize the classification results for each value of $g$. Many of these have been discussed already.

As noted above, Cartan found complete classifications in the cases $g = 1,2,3$.  In particular,
if $g=1$,
then the isoparametric hypersurface $M$ is totally umbilic, and it must be an open subset of a great or small
sphere.  If $g=2$, Cartan \cite{Car2} showed that $M$ must be an open subset of a standard product of two spheres,
\begin{displaymath}
S^k(r) \times S^{n-k-1}(s) \subset S^n, \quad r^2+s^2=1.
\end{displaymath}

In the case $g=3$, Cartan \cite{Car3}
showed that all the principal curvatures must have the same multiplicity
$m=1,2,4$ or 8, and the isoparametric hypersurface must be an open subset of a tube of
constant radius 
over a standard embedding of a projective
plane ${\bf FP}^2$ into $S^{3m+1}$,
where ${\bf F}$ is the division algebra
${\bf R}$, ${\bf C}$, ${\bf H}$ (quaternions),
${\bf O}$ (Cayley numbers), for $m=1,2,4,8,$ respectively.  Thus, up to
congruence, there is only one such family for each value of $m$.

In the case of an 
isoparametric hypersurface with $g=4$ four principal curvatures, 
M\"{u}nzner (see Remark \ref{rem-munzner-thm})  proved that
the principal curvatures can have at most two
distinct multiplicities $m_1,m_2$. 
Next Ferus, Karcher and M\"{u}nzner
\cite{FKM} 
used representations
of Clifford algebras to construct
for any positive integer
$m_1$ an infinite series of isoparametric hypersurfaces
with four principal curvatures having respective multiplicities
$(m_1,m_2)$, where $m_2$ is nondecreasing and
unbounded in each series.  As noted earlier,
these examples are now known as isoparametric hypersurfaces of {\it OT-FKM-type}, because
the construction of Ferus, Karcher and M\"{u}nzner was a generalization of an earlier construction due to Ozeki and Takeuchi \cite{OT}--\cite{OT2}.

Stolz \cite{Stolz}  next proved that the multiplicities $(m_1,m_2)$
of the principal curvatures of an isoparametric
hypersurface with $g=4$ principal curvatures
must be the same as those of the hypersurfaces of
OT-FKM-type, or else $(2,2)$ or $(4,5)$, which are the multiplicities for certain homogeneous examples
that are not of OT-FKM-type  (see the classification 
of homogeneous isoparametric hypersurfaces of Takagi and Takahashi \cite {TT}).

Cecil, Chi and Jensen \cite{CCJ1} then
showed that if the multiplicities of an isoparametric hypersurface with four principal curvatures satisfy the condition
$m_2 \geq 2 m_1 - 1$, then
the hypersurface is of OT-FKM-type.  (A different proof of this result, using isoparametric 
triple systems, was given later by Immervoll \cite{Im}.)

Taken together with known results of Takagi \cite{Takagi} for $m_1 = 1$,
and Ozeki and Takeuchi \cite{OT} for $m_1 = 2$, this result of Cecil, Chi and Jensen handled all
possible pairs of multiplicities except for four cases, the OT-FKM pairs
$(3,4), (6,9)$ and $(7,8)$,and  the homogeneous pair $(4,5)$.

In a series of recent papers, Chi \cite{Chi}--\cite{Chi4} completed the classification of isoparametric hypersurfaces with four principal curvatures.
Specifically, Chi showed that in the cases of multiplicities
$(3,4)$, $(6,9)$ and $(7,8)$, the isoparametric hypersurface must be of OT-FKM-type, and in the case $(4,5)$, it must be homogeneous.

The final conclusion is that an isoparametric hypersurface with $g=4$ principal curvatures must either be of
OT-FKM-type, or else a homogeneous isoparametric hypersurface with multiplicities $(2,2)$ or $(4,5)$.

In the case of an isoparametric hypersurface with $g=6$ principal curvatures, M\"{u}nzner \cite{Mu}--\cite{Mu2} showed
that all of the principal curvatures must have the same multiplicity $m$, and
Abresch \cite{Ab} showed that $m$ must equal 1 or 2.
By the classification of homogeneous isoparametric hypersurfaces due to Takagi and Takahashi \cite{TT},
there is up to congruence only one homogeneous family in each case, $m=1$ or $m=2$.  

These 
homogeneous examples have been shown to be
the only isoparametric hypersurfaces in the case $g=6$ by Dorfmeister and Neher \cite{DN5}
in the case of multiplicity $m=1$, and by Miyaoka \cite{Mi11}--\cite{Mi12} in the case $m=2$
(see also the papers of Siffert \cite{Siffert1}--\cite{Siffert2}, and Miyaoka \cite{Mi25}).

For general surveys on isoparametric hypersurfaces in spheres,
see the papers of Thorbergsson \cite{Th6}, Cecil \cite{Cec9}, and Chi \cite{Chi-survey}.

\noindent Thomas E. Cecil

\noindent Department of Mathematics and Computer Science

\noindent College of the Holy Cross

\noindent Worcester, MA 01610

\noindent email: tcecil@holycross.edu\\

\noindent Patrick J. Ryan

\noindent Department of Mathematics and Statistics

\noindent McMaster University

\noindent Hamilton, Ontario, Canada L8S4K1

\noindent email: ryanpj@mcmaster.ca

\end{document}